\documentclass[12pt]{article}

\usepackage{amsmath}
\usepackage{amssymb}
\usepackage{bm}
\usepackage[dvipdfm]{graphicx}

\def \C {{\mathbb C}}

\def \F {{\mathbb F}}
\def \R {{\mathbb R}}

\def \Z {{\mathbb Z}}
\def \N {{\mathbb N}}

\def \Th {{\rm Th}}

\def \acl {{\rm acl}}

\def \L {{\cal L}}

\def \proof {{\noindent \bf Proof: }}
\def \claim {{\noindent \bf Claim }}
\def \qed 
{
   \hfill \rule {8pt}{8pt} \bigskip
}

\def \F {{\mathbb F}}
\def \C {{\mathbb C}}
\def \R {{\mathbb R}}

\def \Z {{\mathbb Z}}

\def \N {{\mathbb N}}

\def \Th {{\rm Th}}

\def \acl {{\rm acl}}

\def \deg {{\rm deg}}

\def \L  {{\cal L}}

\def \U {{\mathbf U}}
\def \V {{\mathbf V}}
\def \u {{\mathbf u}}
\def \v {{\mathbf v}}
\def \ran {{\rm ran}}

\newtheorem{thm}{Theorem}
\newtheorem{lemma}[thm]{Lemma}
\newtheorem{prop}[thm]{Proposition}
\newtheorem{remark}[thm]{Remark}

\newtheorem{defn}[thm]{Definition}
\newtheorem{cor}[thm]{Corollary}

\title{Notes on a model theory of quantum 2-torus $T_q^2$
for generic $q$}
\author{Masanori ITAI
and Boris ZILBER
}
\date{}

\begin{document}

\maketitle

\abstract{We describe a structure  over the complex numbers
associated with the non-commutative algebra $\mathcal{A}_q$ called quantum 2-tori.
These turn out to have uncountably categorical $L_{\omega_1,\omega}$-theory, and are similar 
to other pseudo-analytic structures considered by the second author. 
The first-order theory of a quantum torus for generic $q$ interprets arithmetic and so is unstable 
and undecidable. But certain interesting reduct of the structure, a quantum line bundle, is  
 superstable. 
}  

\section{Introduction}
Quantum tori have been considered in various contexts;
\begin{enumerate}
\item\label{e1}
Generalizing the definition of
algebraic tori one obtains the notion of {\it quantum $n$-torus} over an abstract field $\F$ as 
the $\F$-algebra ${\cal O}_q((\F^{\times})^n)$ with generators
$x_1^{\pm}, \cdots, x_n^{\pm}$ and the relation
$$x_i x_j = q x_j x_i.$$
See 
 \cite{BG}.

\item \label{e2} More often one is interested in a generalization of  tori in the context of real structure 
on complex manifolds. Then the appropriate generalization to non-commutative setting is based on 
the notion of a $C^*$-algebra. Namely, in the above example we would assume 
that $\F=\C$ and ${\cal O}_q((\F^{\times})^n)$ is a $C^*$-algebra (that is a normed algebra with 
an involution $x\mapsto x^*,$ where $x^*$ is read as an operator adjoint to 
$x$), and the $x_i$ are unitary, that is $x^*=x^{-1}$.  
See \cite{C}  

\item Many more beautiful and important examples can be seen as generalizations of the above, 
including quantum groups.
\end{enumerate}

Classically, in non-commutative geometry one studies representation theory of the algebras in question. 
In case  \ref{e1}, finite-dimensional representations in $\F$-vector spaces, and in case \ref{e2}, 
representations in Hilbert spaces.  One of the main suppositions of non-commutative geometry is that 
unlike commutative case there is no geometric object corresponding to  the quantum algebra, and the best 
we can have in place of the Gel'fand-Naimark duality (co-ordinate algebra -- geometric space) is the correspondence 
between the algebra and its category of representations.

In this paper, however, we 
construct geometric objects, which can be seen as representing the information coded in the algebra of 
a quantum $2$-torus. This is similar to what has been done by the second author in \cite{Z3} in the case 
when $q$ is a root of unity. In that case the appropriate geometric object is a Noetherian Zariski geometry, 
and when the algebra in question is commutative the Zariski geometry is just an affine algebraic variety,  
as in the classical duality. 

The object  $T_q^2(\F)$  constructed in this paper over an algebraically closed field $\F$ hopefully can be 
classified as an analytic Zariski geometry in the sense of \cite{Z2}, but we don't prove this fact here. 
Our main result, apart from construction as such,  is that a
simple $L_{\omega_1,\omega}$-sentence characterizing $T_q^2(\F)$ is categorical in uncountable cardinals. 

\medskip

Let $\Gamma$ denote an infinite cyclic group generated by an element $q\in \F^{\times}$.
We denote ${\cal A}_q$ the non-commutative algebra ${\cal O}_q((\F^{\times})^2)$  with generators written as
$U,U^{-1}, V, V^{-1}$ satisfying
$$
VU = qUV.
$$ 
Our objective is to construct a structure $T_q^2(\F)$ which interprets $U$, $U^{-1}$, $V$, $V^{-1}$ as 
operators acting on it and thus represents
 the algebra ${\cal A}_q.$

In section 2, 
we first construct non-commutative geometric objects
called $\Gamma$-bundles and line-bundles.
Then we construct the {\it quantum} 2-torus $T_q^2(\F)$ 
associated with the algebra ${\cal A}_q$ over 
the algebraically closed field $\F$
having two quantum line-bundles with a pairing function.

After constructing $T_q^2(\F)$, in section 3
we study its properties from model theoretic point of view.

Three main theorems proved in this paper are;
\begin{enumerate}
\item
$T_q^2(\F)$ is axiomatisable by an $L_{\omega_1, \omega}$-sentence which is categorical in uncountable cardinals.
(Theorem \ref{Lomega1})
\item
The first-order theory of line-bundles, which is a reduct of $T_q^2(\F)$, is superstable.
(Theorem \ref{main1})
\item The first-order theory of $T_q^2(\F),$ for $q$ not a root of unity, interprets the ring of integers.
Hence the theory of the quantum 2-torus is undecidable and unstable.
(Theorem \ref{main2})
\end{enumerate}

Prerequisites in model theory is minimal and all found in standard 
text books such as \cite{M} or \cite{TZ}. 
Let $\C, \R, \Z, \N$ denote the set of complex numbers, the set of real numbers,
the set of integers, and the set of natural numbers, respectively.
All the fields considered in this paper are of characteristic zero.

\subsection{Further works}

In a forthcoming work we would like to classify $T_q^2(\F)$ as an analytic Zariski structure,
just like \cite{Z3} classifies a quantum torus at root of unity as Noetherian Zariski structure.
This requires a quantifier elimination statement and a detailed analysis of definability,
this time including $\L_{\omega_1, \omega}$-formulae.

A very important model theoretic next step in the study of quantum tori would be the study of
definable bijections between $T_q^2(\F)$ and $T_{q'}^2(\F)$, analogue of regular isomorphisms
between algebraic varieties in algebraic geometry. Here "definable" assumes definability in
an ambient larger structure over a field $\F$, where the $T_q^2(\F)$ are definable for each
choice of $q$. Such a definable bijection (isomorphisms of tori) at the level of coordinate
algebra ${\cal O}_q((\F^{\times})^n)$ must correspond to a Morita equivalence between algebras.

\section{The quantum 2-torus $T_q^2(\F)$}

Let $\F$ be a field, and $q \in \F^*,$ not a root of unity. 
Consider a $\F$-algebra ${\cal A}_q$ generated by operators
$U, U^{-1}, V, V^{-1}$ satisfying
$$
VU = q UV, \quad UU^{-1} = U^{-1}U = VV^{-1} = V^{-1}V = I.
$$
Let $\Gamma_q = q^\Z = \{ q^n : n \in \Z \}$ be 
a cyclic multiplicative subgroup of $\F^{*}$.
From now on in this section we work in an uncountable 
$\F$-module ${\cal M}$ such that $\dim {\cal M} \geq |\F|$.
Also we drop the subscript $q$ from $\Gamma_q$ when it is clear
from the context.
\subsection{$\Gamma$-sets, $\Gamma$-bundles, line-bundles}

For each pair $(u, v) \in \F^* \times \F^*$,
we will construct two ${\cal A}_q$-modules $M_{|u,v \rangle}$ and $M_{\langle v, u|}$ 
so that
both $M_{|u,v \rangle}$ and $M_{\langle v, u|}$
are sub-modules of ${\cal M}$. 

Before starting the construction of $M_{|u,v \rangle}$ and $M_{\langle v, u|}$
for general $q$, it is important to keep in mind 
the case when $q$ is actually a root of unity.
In such case we can also define $M_{|u,v \rangle}$ and $M_{\langle v, u|}$
in the same manner described in this paper, however $M_{|u,v \rangle}$ and $M_{\langle v, u|}$
will be inter-definable uniformly on the pair $(u, v)$.

We now start the construction with $q$ not a root of unity.
The module $M_{|u, v \rangle}$ is generated by linearly independent
elements labelled 
$\{ \u(\gamma u, v) \in {\cal M}: \gamma \in \Gamma \}$ satisfying 

\begin{equation}\label{UVu}
\begin{array}{ccl}
U & : & \u(\gamma u, v) \mapsto \gamma u \u( \gamma u, v), \\
V & : & \u(\gamma u, v) \mapsto v \u( q^{-1} \gamma u, v), \\
\end{array}
\end{equation}

and also
\begin{equation}\label{UVinvu}
\begin{array}{ccl}
U^{-1} & : & \u(\gamma u, v) \mapsto \gamma^{-1} u^{-1} \u( \gamma u, v), \\
V^{-1} & : & \u(\gamma u, v) \mapsto v^{-1} \u( q \gamma u, v).\\
\end{array}
\end{equation}

Next let $\phi : \F^{*} / \Gamma \rightarrow \F^*$ such that $\phi(x\Gamma) \in x\Gamma$
for each $x\Gamma \in \F^*/\Gamma$.
Put $\Phi$ = $\ran (\phi)$.
We call $\phi$ a {\it choice function} and $\Phi$ the system of representatives.

Set for $\langle u, v \rangle \in \Phi^2$
\begin{equation}
\begin{array}{lcl}
\Gamma \cdot \u(u, v) & := & \{ \gamma \u(u, v) : \gamma \in \Gamma \}, \\
\U_{\langle u, v \rangle} & := & \bigcup_{\gamma \in \Gamma} \Gamma \cdot \u(\gamma u, v)  = 
\{ \gamma_1 \cdot \u( \gamma_2 u, v) : \gamma_1, \gamma_2 \in \Gamma \}.
\end{array}
\end{equation}

And set
\begin{equation}\label{FU}
\begin{array}{lcl}
\U_\phi & := & \bigcup_{\langle u, v \rangle \in \Phi^2} \U_{\langle u, v \rangle} \\
 & = & 
\{ \gamma_1 \cdot \u( \gamma_2 u, v) : 
\langle u, v \rangle \in \Phi^2, \gamma_1. \gamma_2 \in \Gamma \}, \\
\F^* \U_\phi & := &
\{ x \cdot \u(\gamma u, v) : \langle u, v \rangle \in \Phi^2,
x \in \F^*, \gamma \in \Gamma \}. \\
\end{array}
\end{equation}

Note: The notation like $x \cdot \u(\gamma u, v)$ above should be read as
a 4-tuple $(x, \gamma, u, v)$.

We call $\Gamma \cdot \u(u,v)$ a $\Gamma$-{\bf set} over the pair $(u, v)$, 
$\U_\phi$ a $\Gamma$-bundle over  $\F^*\times \F^*/\Gamma$,
and $\F^* \U_\phi$ a {\bf line-bundle} over $\F^*$.
Notice that $\U_\phi$ can also be seen as a bundle inside
$
\bigcup_{\langle u, v \rangle} M_{|u, v \rangle}.
$
Furthermore the line bundle $\F^*\U_\phi$ is closed under the action of 
the operators $U$ and $V$ satisfying the relations (\ref{UVu}) and (\ref{UVinvu}).

We define the module $M_{\langle v, u |}$ generated by 
linearly independent elements labelled 
$\{ \v(\gamma v, u) \in {\cal M} : \gamma \in \Gamma \}$ satisfying 
\begin{equation}\label{UVv}
\begin{array}{ccl}
U & : & \v(\gamma v, u) \mapsto u \v( q\gamma v, u), \\
V & : & \v(\gamma v, u) \mapsto \gamma v \v( \gamma v, u), \\
\end{array}
\end{equation}
and
\begin{equation}\label{UVinvv}
\begin{array}{ccl}
U^{-1} & : & \v(\gamma v, u) \mapsto u^{-1} \v( q^{-1} \gamma v, u), \\
V^{-1} & : & \v(\gamma v, u) \mapsto \gamma^{-1} v^{-1} \v( \gamma v, u). \\
\end{array}
\end{equation}

Similarly a $\Gamma$-set
 $\Gamma \cdot \v(v, u)$ over the pair $(v, u)$, 
a $\Gamma$-bundle $\V_\phi$ over $\F^*/\Gamma\times \F^*$ , 
and $\F^* \V_\phi$ a {\bf line-bundle} over $\F^*$ are defined.

As before the $\Gamma$-bundle $\V_\phi$  can also be seen as a bundle inside
$
\bigcup_{\langle v, u \rangle} M_{\langle v, u |}.
$

In the next section we treat $\F^* \U_\phi$ as an object definable 
in the structure $(\U_\phi, \V_\phi, \F)$. 
For this we introduce an equivalence relation $E$ 
identifying $\gamma \in \Gamma$ as an element of $\F^*$. 
Thus
$\F^* \U_\phi \simeq
(\F^* \times \U_\phi)/E$ where 
for $(x, \gamma_1 \cdot \u( \gamma_3u, v)), 
(x', \gamma_2 \cdot \u( \gamma_3u, v)) \in \F^* \times \U_\phi$ 
define
$$ 
\begin{array}{cc}\label{equiv}
(x, \gamma_1 \cdot \u( \gamma_3u, v)) \sim_E (x', \gamma_2 \cdot \u( \gamma_3u, v)) &
\Longleftrightarrow \\
\exists \gamma \in \Gamma \, ( x' = x \gamma^{-1} 
\wedge \gamma_2 = \gamma \gamma_1 ) &
\end{array}
$$
We will then
consider that the two operators $U$ and $V$ are acting on this definable 
classes $\F^* \U_\phi$ in the next section.

Similarly for $\F^* \V_\phi$ and the actions of $U, V$.

\begin{remark}
Relations (\ref{UVu}) mean that elements $\u(\gamma u, v)$ are 
{\rm eigenvectors}
and $\gamma u$ are 
{\rm eigenvalues} of the operator U.
\end{remark}

Having defined the line bundles $\F^*\U_{\phi}$ and $\F^*\V_{\phi}$, we realize
that any particular properties of the element $q$ or the choice function
$\phi$ are not used. 
This means the following:

\begin{prop}\label{iso1}
Let $\F$, $\F'$ be fields  
and $q \in \F$, $q' \in \F'$ such that
there is an field isomorphism $i$ from $\F$ to $\F'$ sending $q$ to $q'$.
Then $i$ can be extended to an isomorphism from
the $\Gamma$-bundle $\U_{\phi}$ to the $\Gamma'$-bundle $\U_{\phi'}$  
and also from the line-bundle $\F^*\U_{\phi}$ to the line-bundle $(\F^*)'\U_{\phi'}$.
The same is true for the line-bundles $\F^*\V_{\phi}$ and $(\F')^*\V_{\phi'}$.

In particular the isomorphism type of $\Gamma$-bundles and line-bundles does
not depend on the choice function.
\end{prop}

\proof
Let $i$ be an isomorphism from $\F$ to $\F'$ sending $q$ to $q'$.
Set $i (x \cdot \u (\gamma u, v)) = i(x) \cdot \u ( i(\gamma u), i(v))$.
Then this defines an isomorphism from $\F^* \U_\phi$ to $(\F')^* \U_{\phi'}$.
\qed

\subsection{Pairing function and quantum 2-torus}

It is clear from the construction that there is no interactions between 
$\Gamma$-bundles $\U_{\phi}$ and $\V_{\phi}$.
We now introduce the notion of {\it pairing function}
$\langle \cdot \, | \, \cdot \rangle$ 
which plays the r\^{o}le of an {\it inner product}
of two $\Gamma$-bundles $\U_\phi$ and $\V_\phi$:
\begin{equation}\label{pairing}
\langle \cdot \, | \, \cdot \rangle : \Big( \V_\phi \times \U_\phi \Big) 
\cup \Big( \U_\phi \times \V_\phi \Big)
\rightarrow \Gamma.
\end{equation}

We would like two operators $U, V$ to behave like {\it unitary operators} 
on Hilbert space.
This requirement forces us to postulate the following:
\begin{enumerate}
\item $\langle \u(u,v) | \v(v, u) \rangle = 1$,
\item for each $r, s \in \Z$, 
$\langle U^rV^s \u(u, v) | U^rV^s \v(v, u) \rangle = 1$,
\item for $\gamma_1, \gamma_2, \gamma_3, \gamma_4 \in \Gamma$,
$$\langle \gamma_1 \u(\gamma_2u, v) | \gamma_3 \v(\gamma_4v,u) \rangle =
\langle \gamma_3 \v(\gamma_4v, u) | \gamma_1 \u(\gamma_2u,v) \rangle^{-1},$$
\item $\langle \gamma_1 \u(\gamma_2u, v) | \gamma_3 \v(\gamma_4v,u) \rangle =
\gamma_1^{-1} \gamma_3 \langle \u(\gamma_2u, v) | \v(\gamma_4v,u) \rangle$, and
\item for $v' \not\in \Gamma \cdot v$ or $u' \not\in \Gamma \cdot u$,
$
\langle q^s \v(v',u) | q^r \u(u',v) \rangle 
$
is not defined. 
\end{enumerate}

\begin{prop}
The {\bf pairing} function (\ref{pairing}) defined above satisfies 
the following:
for any $m, k, r, s \in \N$ we have
\begin{equation}\label{pairing1}
\langle q^s \v(q^mv,u) | q^r \u(q^ku,v) \rangle =
q^{r-s-km}
\end{equation}
and
\begin{equation}\label{pairing2}
\langle q^r \u(q^ku,v) | q^s \v(q^mv,u) \rangle =
q^{km+s-r} =
\langle q^s \v(q^mv,u) | q^r \u(q^ku,v) \rangle^{-1}.
\end{equation}
\end{prop}

\proof
We only prove (\ref{pairing1}). For this, it is enough to notice that 
for each $r, s \in \Z$, 
$\langle U^rV^s \u(u, v) | U^rV^s \v(v, u) \rangle = 1$ implies that
$$\langle \u(q^ru, v) | \v(q^sv, u) \rangle = q^{rs}.$$ 
\qed

Finally we are ready to define the notion of quantum 2-torus:

\begin{defn}\label{q-torus}
We call the multi-sorted structure 
$(\U_\phi,  \V_\phi, \langle \cdot \, | \, \cdot \rangle, \F )$ with actions
$U$, $V$ with $U$ satisfying (\ref{UVu}),  
$V$ satisfying (\ref{UVv}) 
and $\langle \cdot \, | \, \cdot \rangle$ 
a pairing function defined as above a {\bf quantum 2-torus} $T_q^2(\F)$ over the 
field $\F$.
\end{defn}

From Proposition \ref{iso1} we know that the structure of the line-bundles does
not depend on the choice function. The next proposition tells us that
the structure of the quantum 2-torus $T_q^2(\F)$ depends only on $\F$, $q$
 and not on the choice function.

\begin{prop}[cf. Proposition 4.4, \cite{Z1}]
Given $q\in \F^*$ not a root of unity, 
any two structures of the form $T_q^2(\F)$ are isomorphic
over $\F$. In other words, the isomorphism type of $T_q^2(\F)$ does not
depend on the system of representatives $\Phi$.
\end{prop}

\proof
Let $\phi, \psi$ be two choice functions
of $\F^*/\Gamma$. Consider two structures 
$(\U, \V)_\phi = (\U_\phi, \V_\phi)$ and 
$(\U, \V)_\psi = (\U_\psi, \V_\psi)$.
We show that these two structures are isomorphic over $\F$.

Suppose $\phi$ picks $\langle u_g, v_g \rangle$ from $\F^*/\Gamma$ and
$\psi$ picks $\langle u_0, v_0 \rangle$ from the same coset of $\langle u_g, v_g \rangle$. 

Consider the bases 
$\{ \u(q^ku_g, v_g) : k \in \Z \}$ of $U_{\langle u_g, v_g \rangle}$ and 
$\{ \v(q^kv_g, u_g) : k \in \Z \}$ of $V_{\langle u_g, v_g \rangle}$ in the
structure $(\U_\phi, \V_\phi)$.

Since $\langle u_0, v_0 \rangle$ and $\langle u_g, v_g \rangle$ are in the same coset
of $\F^*/\Gamma$ there are $s, t \in \Z$ such that $u_0 = q^su_g, v_0 = q^tv_g$. 

We now want to transfer the structure of $\U_{\langle u_g, v_g \rangle}$ and 
$\V_{\langle u_g, v_g \rangle}$
to $\U_{\langle u_0, v_0 \rangle}$ and $\V_{\langle u_0, v_ \rangle}$ 
respectively as follows. Set
\begin{itemize}
\item $\u(u_0, v_0) := q^{st}\u(q^su_g, v_g)$,
\item $\u(q^ku_0, v_0) := v_0^kV^{-k}\u(u_0, v_0)$,
\item $\v(v_0, u_0) := \v(q^tv_g, u_g)$,
\item $\v(q^kv_0, u_0) := u_0^{-k} U^k \v(v_0, u_0)$,
\end{itemize}
where $k \in \Z$.
First notice that we have
$$
\begin{array}{lcl}
\u(q^ku_0,v_0) & = & v_0^k V^{-k} \u (u_0, v_0) \\
  & = & (q^t v_g)^k V^{-k} (q^{st} \u (q^s u_g, v_g)) \\
  & = & q^{kt} v_g^k q^{st} v_g^{-k} \u(q^{s+k}u_g, v_g) \\
  & = & q^{kt + st} \u(q^{k+s}u_g, v_g), 
\end{array}
$$
and
$$
\begin{array}{lcl}
\v(q^kv_0,u_0) & = & u_0^{-k} U^k \v(v_0, u_0) \\
  & = & q^{-ks}u_g^{-k} U^k \v( g^tv_g, u_g) \\
  & = & q^{-ks}u_g^{-k} u_g^k \v( g^{t+k}v_g, u_g) \\ 
  & = & q^{-sk}\v(q^{k+t}v_g, u_g). \\
\end{array}
$$
From these relations we see that
the operators $U$ and $V$ act on the set 
$\{ \u(q^ku_0, v_0) : k \in \Z \}$ according to the definition of
$\U_{\phi}$, that is
$$\begin{array}{lcl}
U(\u(u_0, v_0)) & = & U(q^{st}\u(q^su_g, v_g)) \\
 & = & q^{st} U(\u(q^su_g, v_g)) \\
 & = & q^{st}q^su_g\u(q^su_g,v_g) \\
 & = & u_0 q^{st}\u(q^su_g,v_g) \\
 & = & u_0 \u(u_0,v_0), \\
V(\u(u_0, v_0)) & = & V(q^{st}\u(q^su_g, v_g)) \\
 & = & v_0 \u(q^{-1}u_0,v_0), \\
U(\u(q^ku_0,v_0)) & = & U(v_0^kV^{-k}\u(u_0,v_0)) \\
  & = & q^ku_0 \u(q^ku_0, v_0), \\
V(\u(q^ku_0,v_0)) & = & V( q^{kt + st} \u(q^{k+s}u_g, v_g)) \\
  & = & v_0 \u(q^{k-1}u_0, v_0). \\ 
\end{array}$$

By similar calculations we see that the operators $U$ and $V$ act on the set 
$\{ \v(q^kv_0, u_0) : k \in \Z \}$ according to the definition of $\V_\phi$.

Finally from the following relations we see that we can properly transfer
the pairing function from $(\U_{\langle u_g, v_g \rangle}, \V_{\langle v_g, u_g \rangle})$
to $(\U_{\langle u_0, v_0 \rangle}, \V_{\langle v_0, u_0 \rangle})$:
$$\langle \v(v_0,u_0)|\u(u_0,v_0)\rangle =
\langle \v(q^tv_g, u_g)|q^{st}\u(q^su_g, v_g) \rangle = q^{st-st}=1$$
and
$$\begin{array}{lcl}
\langle \v(q^mv_0,u_0)|\u(q^ku_0,v_0) \rangle & = & 
\langle q^{-sm}\v(q^{m+t}v_g,u_g)|q^{st+kt}\u(q^{k+s}u_g,v_g) \rangle \\
 & = & q^{st+kt-(-sm)-(m+t)(k+s)} \\
 & = & q^{-mk}.
\end{array}$$
We have now shown that the two structures 
$(\U_{\langle u_g, v_g \rangle}, \V_{\langle v_g, u_g \rangle})$
and $(\U_{\langle u_0, v_0 \rangle}, \V_{\langle v_0, u_0 \rangle})$ are isomorphic.
Therefore so are the two structures $(\U, \V)_\phi$ and $(\U, \V)_\psi$.
\qed

From the above proposition we have as a corollary to Proposition \ref{iso1}
the following:

\begin{cor}\label{isomorph}
Let $\F$, $\F'$ be fields  
and $q \in \F$, $q' \in \F'$ such that
there is a field isomorphism $i$ from $\F$ to $\F'$ sending $q$ to $q'$.
Then $i$ can be extended to an isomorphism 
from the quantum 2-torus
$T_q^2(\F)$ to the quantum 2-torus $T_{q'}^2(\F')$.

In particular the isomorphism type of quantum 2-torus does
not depend on the choice function.

\end{cor}

From now on we drop the subscript $\phi$ from line-bundles $\U_\phi$, 
$\V_\phi$ and write simply as 
$T_q^2(\F) = (\U, \V, \langle \cdot \, | \, \cdot \rangle, \F)$.

%
%
%
%

\section{The model theory of quantum 2-torus over algebraically closed field}

We now study the model theory of quantum 2-tori $T_q^2(\F)$ under the assumption
that $\F$ is algebraically closed.

After introducing an appropriate language, we shall prove three theorems in this
section;
\begin{enumerate}
\item
$T_q^2(\F)$ is axiomatisable by an $L_{\omega_1, \omega}$-sentence $\Psi$ which is categorical 
in uncountable cardinals.
(Theorem \ref{Lomega1})
\item
The first-order theory of line-bundles, which is a reduct of $T_q^2(\F)$, is superstable.
(Theorem \ref{main1})
\item The first-order theory of $T_q^2(\F),$ for $q$ not a root of unity, 
interprets the ring of integers.
Hence the theory of the quantum 2-torus  is undecidable and unstable.
(Theorem \ref{main2})
\end{enumerate}

\subsection{The language for quantum 2-tori}

To define the sentence $\Psi$ we introduce a language $\L_q$ which is the
language for multi-sorted structure 
$(\U, \V, \langle \cdot \, | \, \cdot \rangle, \F)$ ; \\
$$\L_q = \{ +, \cdot, \U, \V, U, V, \langle \cdot \, | \, \cdot \rangle, \F, \Gamma, \pi \}$$
where $+, \cdot$ defined
on $\F$ and $\Gamma \subset \F$. Furthermore $U, V$ are operators acting
on $\U$ and $\V$. Each $\gamma \in \Gamma$ acts on $\U$ and $\V$. 
Also $\pi$ is a function symbol which will be interpreted as a surjection
from $\U$ onto $\F^* \times \F^*/\Gamma$ and from $\V$ onto $\F^* \times \F^*/\Gamma$.

\subsection{$T_q^2(\F)$ is $L_{\omega_1, \omega}$-categorical 
in uncountable cardinals}

Here we define the $\L_{\omega_1, \omega}$-sentence $\Psi$ in $\L_q$
describing the quantum 2-torus $T_q^2(\F) = (\U, \V, \langle \cdot \, | \, \cdot \rangle, \F)$.
Then we show that the sentence $\Psi$ is categorical in uncountable cardinals.

Recall how we treat $\F^* \U$ as an object definable in the structure 
$(\U, \V, \langle \cdot \, | \, \cdot \rangle, \F)$. 
Similarly for $\F^* \V$.
We extend that the two operators $U$ and $V$ are acting on 
these definable sets $\F^* \U$ and $\F^* \V$.

Let $\Psi$ be the $L_{\omega_1, \omega}$-sentence stating that
\begin{enumerate}
\item $\F$ is an algebraically closed field of characteristic zero,
\item $q \in \F$ and not a root of unity,
\item $\Gamma$ is a multiplicative subgroup of $\F$ generated by $q$, i.e., $\Gamma \simeq q^{\Z}$,
\item $\pi$ is surjective from $\U$ onto $\F^* \times \F^*/\Gamma$,
\item for each $\gamma \in \Gamma, (u, v) \in \F^* \times \F^*/\Gamma$, 
$\pi^{-1}(\gamma u,v) \subset \U$ is generated by an element and
 \begin{itemize}
  \item for each $\u \in \pi^{-1}(\gamma u,v)$ and $\gamma' \in \Gamma, \gamma' \u \in \pi^{-1}(\gamma u,v)$,
  \item for each $\u \in \pi^{-1}(\gamma u,v)$ and $x \in \F^*, x\u \in \F^* \pi^{-1}(\gamma u,v)$,
 \end{itemize}
\item $\F^* \U$, $\F^*V$ are $\F$-modules,
\item operators $U, V$ act on $\F^*\U$ and $\F^*\V$ according to (\ref{UVu}) and (\ref{UVv}), 
more precisely,
 \begin{itemize}
  \item for each $\u \in \pi^{-1}(\gamma u,v), x \in \F^*$ we have 
        $U(x, \u) \in \F^*U$, $\gamma u \u \in \pi^{-1}(\gamma u,v)$ and
        $U(x, \u) = \gamma u \u$,
  \item for each $\u \in \pi^{-1}(\gamma u,v), x \in \F^*$ we have 
        $V(x, \u) \in \F^*\V$,  there exists 
        $\u' \in \pi^{-1}(q^{-1}\gamma u, v)$ and $V(x, \u) = xv\u'$,
 \end{itemize}
\item the properties of the pairing function, more precisely,
 \begin{itemize}
  \item for any $\u \in \pi^{-1}(u,v), \v \in \pi^{-1}(v,u)$, 
   \begin{itemize}
    \item $\langle \u,\v \rangle = 1$,
    \item for any $r, s \in \Z$, 
    $\langle U^rV^s(\u), U^rV^s(\v) \rangle = 1$, \\
    where $1$ is the multiplicative identity element of $\F$.
   \end{itemize}
  \item for any $\gamma_1, \gamma_2, \gamma_3, \gamma_4 \in \Gamma$ and 
  $\u \in \pi^{-1}(\gamma_2, v), \v \in \pi^{-1}(\gamma_4 v, u)$,
  \begin{itemize}
   \item $\langle \gamma_1 \u, \gamma_3 \v \rangle =
   \langle \gamma_3 \v, \gamma_1 \u \rangle^{-1}$ and,
   \item $\langle \gamma_1 \u, \gamma_3 \v \rangle = \gamma_1^{-1} \gamma_3 
   \langle \u, \v \rangle$. 
  \end{itemize}
 \end{itemize}
\end{enumerate}

\begin{lemma}
Let $M = (\U, \V, \langle \cdot \, | \, \cdot \rangle, \F)$ be an infinite 
$\L_q$-structure satisfying $\Psi$. 
Then $M$ is a quantum 2-torus over $\F$.
\end{lemma}

\proof
Let $\F$ be an infinite algebraically closed field of characteristic zero,
$q \in \F$ and $\Gamma = q^\Z \subset \F$ given by $\Psi$.
Then $\U$ and $\V$ are the $\Gamma$-bundles defined over $\F$ and $q$ where
operators $U, V$ satisfy relations (\ref{UVu}), (\ref{UVinvu}), (\ref{UVv}) and (\ref{UVinvv}).
We then construct line-bundles $\F^*\U$ and $\F^*\V$. Since the sentence $\Psi$ describes the
properties of the pairing function, from Proposition \ref{iso1} we
see that $M$ is a quantum 2-torus over $\F$.
\qed

Thus Corollary \ref{isomorph} gives us

\begin{thm}\label{Lomega1}
The $L_{\omega_1, \omega}$-sentence $\Psi$ is categorical in uncountable cardinals.
\end{thm}

\proof
Notice that any uncountable model of $\Psi$ is a quantum 2-torus $T_q^2(\F)$
over an uncountable algebraically closed field $\F$ of characteristic zero.
Since all such uncountable fields $\F$ are isomorphic once we fix the cardinality, 
so are the quantum 2-tori $T_q^2(\F)$ over such fields $\F$.
Hence $\Psi$ is categorical in uncountable cardinals.
\qed

\subsection{Superstability of $\Th(\F, +, \cdot, 0, 1, \Gamma)$}

From now on we study the first-order theoretic properties of quantum 2-tori.
In this subsection we show that the first-order 
theory of $(\F, +, \cdot, 0, 1, \Gamma)$ is axiomatizable and superstable. 
Key idea is that the predicate
$\Gamma(x)$ describes the property of the set $q^\Z$ as a multiplicative
subgroup with the following Lang-type property. 

\begin{defn}[Definition 2.3 \cite{P}]\label{LTD}
Let $K$ be an algebraically closed field, and $A$ a commutative algebraic
group over $K$ and $\Gamma$ a subgroup of $A$. 
We say that $(K, A, \Gamma)$ is of {\rm Lang-type} if
for every $n < \omega$ and every subvariety $X$ (over $K$)
of $A^n = A \times \cdots \times A$ ($n$ times),
$X \cap \Gamma^n$ is a finite union of cosets of 
subgroups of $\Gamma^n$.
\end{defn}

The Lang-type property gives us 
\begin{prop}[Proposition 2.6 \cite{P}]\label{LTP}
Let $K$ be an algebraically closed field, $A$ a commutative algebraic group
over $K$, and $\Gamma$ a subgroup of $A$.
Then $(K, A, \Gamma)$ is of Lang-type if and only if 
$\Th(K, +, \cdot, \Gamma, a)_{a \in K}$ is stable and $\Gamma(x)$ is
one-based.
\end{prop}
Here $\Gamma(x)$ is {\it one based} means that for every $n$ and
every definable subset $X \subset \Gamma^n$, $X$ is a finite
boolean combination of cosets of definable subgroups of
$\Gamma^n$.

With the above Definition \ref{LTD} and Proposition \ref{LTP} in mind, we axiomatize 
the properties of $(\F, +, \cdot, \Gamma )$ as follows;

\bigskip
\noindent
{\bf Axioms for $(\F, +, \cdot, \Gamma )$}
\begin{itemize}
\item[A. 1]
$\Gamma$ satisfies the first order theory of a cyclic group with generator $q,$ 
\item[A. 2] (Lang-type) for every $n$ and every variety $X$ of $(\F^*)^n$, 
$X \cap \Gamma^n$ is a finite union of cosets of definable subgroups
of $\Gamma^n$.
\end{itemize}

Let $T_{\F, \Gamma}$ denote the set of all logical consequences of the 
axioms for $\Gamma$ and ${\rm ACF}_0$ axioms for the algebraically closed
fields of characteristic zero.

\begin{lemma}\label{Lang}
The {\rm Lang-type} property of $(\F, +, \cdot, \Gamma)$ is witnessed by its 
first-order theory.
\end{lemma}

\proof
We may suppose $X$ is irreducible.
Each such variety $X \subset (\F^*)^n$
is definable by an irreducible polynomial $f(x_1, \cdots, x_n)$
over $\F^*$. Definable cosets of $\Gamma^n$ are of the form 
$\overline{\gamma} \Gamma^n = 
\gamma_1 \Gamma \times \cdots \times \gamma_n \Gamma$
where $\gamma_1, \cdots, \gamma_n \in \Gamma(\F)$. 
Hence the sentence "$X \cap \Gamma^n$ {\it is a finite union of cosets of 
definable subgroups}" is expressed as
$$
\begin{array}{ll}
(f(x_1, \cdots, x_n) = 0) \wedge \Gamma(x_1) \wedge \cdots \wedge
\Gamma(x_n) & \longleftrightarrow \\ 
\bigvee_{i=1}^{N_f} \varphi_i(x_1, \cdots, x_n). & \\
\end{array}  
$$
Where each $\varphi_i(x_1, \cdots, x_n)$ defines a coset.
Crucial point here is that the number $N_f$ of the bound of cosets is 
computable for each polynomial $f$. 
For this note first that for any $k \in \N$
the number of cosets of $q^{k\Z}$ in $q^\Z$ is $k$.
Suppose 
$$
f(x_1, \cdots, x_n) = \sum_{i=0}^{\deg(f)} \overline{a}_i 
\overline{x}_i^{m_1} ,
$$ 
where each $m_i$ is a multi index. Let $M_i$ be the sum of
multi index $m_i$. Then the bound $N_f$ of number of cosets
is $\deg(f) \cdot \sum_{i=0}^{\deg(f)} M_i$.
Therefore the Lang-type property is first-order.
\qed

\begin{prop}\label{CompFGamma}
$T_{\F, \Gamma}$ is complete. Hence $T_{\F, \Gamma} =
\Th(\F, +, \cdot, \Gamma)$.
\end{prop}

\proof
Consider a saturated model $(\F, +, \cdot, \Gamma, q)$ 
of $T_{\F, \Gamma}$.
Set $\Gamma(\F) = \{ x \in \F :  \F \models \Gamma(x) \}$. 
Let $q$ be an element of $\F$ interpreting the constant. 
By the axioms
for $\Gamma$, $q^\Z \subset \Gamma(\F) \subset \F$ .

Consider a complete type $t_0(x)$ generated by the following set of formulas, 
$$
t(x) = \{ \Gamma(x), \, \exists y (x = qy), \, \exists y ( x = q^2y), \, \cdots \}.
$$
By saturation there exists $\gamma_0 \in \Gamma(\F)$ realizing $t_0(x)$.  
Clearly, $\gamma_0 \notin q^\Z$.
Suppose elements $\gamma_0, \cdots , \gamma_i \in \Gamma(\F)$ 
have been defined.
Let $t_{i+1}(x)$ be a complete type generated by the type $t(x)$ and the set
$$
\{ x \neq \gamma_0^{n_0} , \cdots,  x \neq \gamma_i^{n_i}  
\, : \, n_0, \cdots, n_i \in \Z \}.
$$
From saturation, we have $\gamma_{i+1} \in \Gamma(\F)$ such that 
$$
\gamma_{i+1} \notin \bigcup_{l = 0}^{i} \gamma_l^\Z.   
$$
In this way by saturation as before 
we see that there exist 
$\gamma_0, \gamma_1, \cdots,  \gamma_i , $ 
$\cdots \in \Gamma(\F)$ $(i < |\F|)$
such that 
$$\Gamma(\F) = q^\Z \cup \bigcup_{i < |\F|}  \gamma_i^\Z.$$
Now take two saturated models $(\F, +, \cdot, \Gamma, q)$ and $(\F', +, \cdot, \Gamma', q')$ 
of $T_{\F, \Gamma}$ of the same cardinality.  There is an isomorphism $i$ from $\F$ to $\F'$
sending $q$ to $q'$.
By the above formula for $\Gamma(\F)$ and the back-and-forth argument 
we can extend $i$ to have that $\Gamma(\F) \simeq \Gamma'(\F')$. Hence 
$(\F, +, \cdot, \Gamma, q)$ and $(\F', +, \cdot, \Gamma', q')$
are isomorphic as saturated models of $T_{\F, \Gamma}$. This completes the
proof of the completeness of the theory $T_{\F, \Gamma}$. 
\qed

\begin{thm}\label{main}
$T_{\F, \Gamma}$ is superstable.
\end{thm}

\proof
Notice first that the multiplication of $\F$ is an algebraic group and
$(\F, \cdot, \Gamma)$ is of the Lang-type by A. 2 above. 
Thus by Proposition 2.6 of \cite{P}, we see that $T_{\F, \Gamma}$ is
at least stable. $T_{\F, \Gamma}$ is in fact superstable since 
\begin{enumerate}
\item 
the stability spectrum of $T_{\F, \Gamma}$ is the same as
that of $T_{\Gamma(\F)}$, the theory of restriction of $(\F, +, \cdot, \Gamma)$ to $\Gamma(\F)$. 
Let $C \subset \F$. Observe first that there is only one complete 1-type over $C$ in $T_{\F, \Gamma}$,
which is realized by elements in $\F - \acl_{\rm F}(\Gamma(\F) \cup C)$ where $\acl_{\rm F}$ 
is the field-theoretic algebraic closure.  
Hence the cardinality of complete 1-types in $T_{\F, \Gamma}$ is bounded by the
cardinality of the complete 1-types in $T_{\Gamma(\F)}$. Thus they have the same stability spectrum.
\item $T_{\Gamma(\F)}$ is superstable. 
For $q$ transcendental, this is Theorem 1 in Section 5 of \cite{Z7}.
Combined with Proposition \ref{LTP}, it is easy to extend it to arbitrary $q$.
\end{enumerate}
\qed

\subsection{Superstability of the line-bundle}

In this subsection we show that the first-order theory of the line-bundle $(\U, \F)$ is superstable.

Recall the $\L_{\omega_1, \omega}$-sentence $\Psi$ describing the properties of
the quantum 2-torus $T_q^2(\F)$. We now investigate the sentence $\Psi$ from
the first-order theoretic point of view.

Let $\L_q' = \L_q - \{ \langle \cdot \, | \, \cdot \rangle \}$.
Let $\Th(\U, \F)$ denote the first-order $\L_q'$-theory of the line bundle $(\U, \F)$.
Unlike $\L_{\omega_1, \omega}$-sentence $\Psi$, in $\Th(\U, \F)$ we can only say
that $\Gamma$ satisfies the first-order theory of a cyclic group with generator $q$. 

Let $M$ be a model of $\Th(\U, \F)$,  then we have;
for each $\gamma \in \Gamma$ and $(u, v) \in \F^* \times \F^*/\Gamma$,
$\pi^{-1}(\gamma u, v)$, which is denoted as $\U_{(\gamma u,v)}$, is a subset of 
the $\Gamma$-bundle $\U$ that is generated by an element.
We also have $\U_{(\gamma u,v)} \subsetneq \F^* \U_{(\gamma u,v)}$
wher $\F^* \U_{(\gamma u,v)}$ is a subset of the line-bundle $\F^* \U$.
Furthermore we have
$$
\U = \bigcup_{\gamma, (u,v)} \U_{(\gamma u,v)}, \quad 
\F^* \U = \bigcup_{\gamma, (u,v)} \F^* \U_{(\gamma u,v)}.
$$

\begin{prop} \label{TUF_is_superstable}
The first-order theory of line-bundle $(\U, \F)$ is superstable.
\end{prop}

\proof
We show that the first-order theory of line-bundle $(\U, \F)$ 
is superstable in two steps; 

1) show that the theory $\Th(\U, \F)$ of line-bundle 
$(\U, \F)$ is {\it prime} over the theory $T_{\F, \Gamma}$;
i.e., any isomorphism between two models of
$T_{\F, \Gamma}$ can be extended to an
isomorphism between two models of $\Th(\U, \F)$,

2) since the theory $T_{\F, \Gamma}$ is superstable
 (Theorem \ref{main})
the theory $\Th(\U, \F)$ of line-bundle
$(\U, \F)$ is superstable as well.

{\bf Proof} of 1): Since the theory $T_{\F, \Gamma}$ is complete
 (Proposition \ref{CompFGamma}), 
we may assume that
$(\F_1, \Gamma) = (\F_2, \Gamma)$.
We show that two saturated models $M_1$ and $M_2$ of $\Th(\U, \F)$
with the same cardinality are isomorphic as line-bundles over
$\F^* \times \F^*/\Gamma$.

Let $\U_1$ and $\U_2$ denote the $\Gamma$-bundles over 
$\F^* \times \F^*/\Gamma$ in $M_1$ and $M_2$ respectively. 
Similarly, let $\F^* \U_1$ and
$\F^* \U_2$ denote the line-bundles over $\F^* \times \F^*/\Gamma$ in 
$M_1$ and $M_2$ respectively.

Take $(u, v), (u', v') \in \F^* \times \F^*/\Gamma$.  
Consider vectors $\u(u, v) \in M_1$ and $\u(u', v') \in M_2$ respectively. 
The $\Gamma$-sets $\Gamma \cdot \u(u,v)$ in $M_1$ and $\Gamma \cdot \u(u', v')$ in $M_2$ are
isomorphic.  Therefore those 1-dimensional submodules generated by
$\Gamma \cdot \u(u, v)$ and $\Gamma \cdot \u(u', v')$ are isomorphic as well.
Hence we see that  
$\F^* \U_{(u,v)}$ in $M_1$ and $\F^* \U_{(u', v')}$ in $M_2$ are isomorphic.

Now by applying the operator $V$ to $\u(u, v)$ in $M_1$ and $\u(u', v')$ in $M_2$
we move to other $\Gamma$-sets $\Gamma \cdot V(\u(u, v))$ in 
$M_1$ and  $\Gamma \cdot V(\u(u', v'))$ in $M_2$ respectively. 
Then two 1-dimensional submodules generated by
$\Gamma \cdot V(\u(u, v))$ and $\Gamma \cdot V(\u(u', v'))$ are isomorphic as well.

In this way we see that the $\Gamma$-bundle in $M_1$
and the $\Gamma$-bundle in $M_2$ are isomorphic. 
Then we extend this isomorphism to an isomorphism between
the line-bundles $M_1$ and $M_2$. This completes the proof of primeness.

{\bf Proof} of 2):  By the primeness shown in 1) 
we see that any realization of a type in $\Th(\U, \F)$ is
fixed by automorphisms of models of $T_{\F, \Gamma}$, hence 
the cardinality of types in
$\Th(\U, \F)$ is bounded by the cardinality of
the types in $T_{\F, \Gamma}$. Therefore by Theorem \ref{main},
the theory $\Th(\U, \F)$ of line-bundle $(\U, \F)$ is superstable.
\qed

It follows immediately that we have the following second main theorem;

\begin{thm}\label{main1}
\begin{enumerate}
\item The first-order theory of line-bundle $(\V, \F)$ is superstable,
\item The first-order theory of line-bundles $(\U, \V, \F)$ is superstable.
\end{enumerate}
\end{thm}

\subsection{Arithmetic in the theory of quantum 2-torus}

In this subsection we show that with the pairing function the ring of integers
can be defined in $\Gamma$. In this regard it is similar to the theory of pseudo-exponentiation, 
the model theory of which can successfully be investigated ``modulo  arithmetic'' as in \cite{Z5} or \cite{KZ}. 

First we may identify $(\Gamma, \cdot)$ with $(\Z, +)$ via the correspondence
$$q^r \mapsto r.$$ 
This gives us immediately a definable addition $+$ on $\Z$ by the exponential law.

A definable multiplication $\times$ on $\Gamma$ is defined as follows with the pairing function.
Fix $u, v \in \F$, $\u(u, v) \in \U$, and $\v(v, u) \in \V$ which satisfy
$$
\langle \u(u, v) | \v(v, u) \rangle = 1.
$$
Then given $\alpha, \beta, \gamma \in \Z$, set by (\ref{pairing2})
$$
\alpha \times \beta = \gamma \quad \text{if and only if} \quad
\langle q^\alpha \u(u, v) | q^\beta \v(v, u) \rangle = q^{\gamma}.
$$

Let $\oplus$ and $\otimes$ be the pull-backs of $+$ and $\times$ in $(\Gamma, \cdot, 1)$ respectively
via the above correspondence. We then have 
\begin{prop}
The two operations $\oplus$ and $\otimes$ defined above are commutative,
i.e., for $\gamma_1, \gamma_2 \in \Gamma$ we have
\begin{enumerate}
\item $\gamma_1 \oplus \gamma_2 = \gamma_2 \oplus \gamma_1$,
\item $\gamma_1 \otimes \gamma_2 = \gamma_2 \otimes \gamma_1$.
\end{enumerate}
\end{prop}

\begin{thm}\label{main2} (i)  With the pairing function,
within $(\Gamma, \cdot, 1, q)$ we can define $(\Gamma, \oplus, \otimes, 1, q)$
and $(\Gamma, \oplus, \otimes, 1, q) \simeq (\Z, +, \cdot, 0, 1)$.
Hence the theory of the quantum 2-torus 
$(\U,  \V, \langle \cdot \, | \, \cdot \rangle, \F, \Gamma)$
is undecidable and unstable.

(ii) The non-elementary theory of the quantum 2-torus
$(\U,  \V, \langle \cdot \, | \, \cdot \rangle, \F, \Gamma)$
over fixed $\Gamma=q^{\Z}$ is categorical in uncountable cardinalities.
\end{thm}

\proof (i) This is essentially the Proposition above.

\noindent
(ii) This is a direct corollary of the two statements

\claim 1. The  non-elementary theory of $(\F,\Gamma),$ the field with a distinguished fixed subgroup, 
is categorical in uncountable cardinalities.

\claim 2. The quantum 2-torus $(\U, \V, \F, \Gamma)$
with the pairing function is prime over  $(\F,\Gamma).$

For Claim 1, note first that the first-order theory of $\F$ is uncountable categorical
since $\F$ is algebraically closed.  For $\Gamma \simeq q^\Z$, what we cannot
state in the first-order theory is 
$x \in \Gamma \longleftrightarrow \exists n \in \Z  \, (x = q^n)$. This is expressible
in the non-elementary theory. Thus the non-elementary theory of $(\F, \Gamma)$
is categorical in uncountable cardinalities.

Claim 2 is in fact part of the proof of Proposition~\ref{TUF_is_superstable}. 
\qed 



DEPARTMENT OF MATHEMATICAL SCIENCES,
TOKAI UNIVERSITY, HIRATSUKA, JAPAN

{\it E-mail address}: itai@tokai-u.jp

\bigskip

MATHEMATICAL INSTITUTE, OXFORD UNIVERSITY, UNITED KINGDOM

{\it E-mail address}: Boris.Zilber@maths.ox.ac.uk
\end{document}